\documentclass[article,onefignum,onetabnum]{siamonline220329}

\usepackage{amsfonts}
\usepackage{graphicx}
\usepackage{epstopdf}
\ifpdf
  \DeclareGraphicsExtensions{.eps,.pdf,.png,.jpg}
\else
  \DeclareGraphicsExtensions{.eps}
\fi

\usepackage{enumitem}
\setlist[enumerate]{leftmargin=.5in}
\setlist[itemize]{leftmargin=.5in}

\newsiamthm{remark}{Remark}
\newtheorem{thm}{Theorem} 
\newtheorem{assumption}{Assumption} 
\crefname{assumption}{Assumption}{Assumptions}

\headers{THE ETKF WITH MULTIPLICATIVE COVARIANCE INFLATION}{K. TAKEDA AND T. SAKAJO}

\title{Uniform error bounds of the ensemble transform Kalman filter for infinite-dimensional dynamics with multiplicative covariance inflation\thanks{Submitted to the editors February 6, 2024.
\funding{The first author is supported by RIKEN Junior Research Associate Program and JST SPRING JPMJSP2110. The second author is supported by JST MIRAI JPMJMI22G1.}}}

\author{Kota Takeda\thanks{Department of Mathematics, Kyoto University, Kitashirakawa Oiwakecho, Sakyo-ku, Kyoto, Kyoto, 606--8502, Japan, and RIKEN Center for Computational Science, 7-1-26, Minatojima Minamimachi, Chuo-ku, Kobe, Hyogo, 650--0047, Japan (\email{takeda.kota.53r@st.kyoto-u.ac.jp}).}
\and Takashi Sakajo\thanks{Department of Mathematics, Kyoto University, Kitashirakawa Oiwakecho, Sakyo-ku, Kyoto, Kyoto, 606--8502, Japan (\email{sakajo@math.kyoto-u.ac.jp}).}}

\usepackage{amsmath, mathrsfs, amsfonts, mathtools, amssymb, amscd, ascmac}
\usepackage{bm}

\usepackage{algorithm,algorithmic}

\newcommand{\E}{\mathbb{E}}
\newcommand{\F}{\mathcal{F}}
\renewcommand{\H}{\mathcal{H}}

\renewcommand{\L}{\mathcal{L}}

\newcommand{\N}{\mathbb{N}}
\newcommand{\Prob}{\mathbb{P}}
\newcommand{\R}{\mathbb{R}}

\newcommand{\U}{\mathcal{U}}

\newcommand{\Y}{\mathcal{Y}}

\newcommand{\enorm}[1]{|#1|_2}
\newcommand{\opnorm}[1]{|#1|_\L}
\newcommand{\wc}{{}\cdot{}}
\newcommand{\bracket}[2]{\left\langle #1, #2 \right\rangle}
\newcommand{\trace}[1]{\operatorname{Tr}\left(#1\right)}
\DeclareMathOperator{\ran}{Ran}
\DeclareMathOperator{\cov}{Cov}

\definecolor{Orange}{HTML}{EB8C00}

\ifpdf
\hypersetup{
  pdftitle={Uniform error bounds of the ensemble transform Kalman ﬁlter for1
  inﬁnite-dimensional dynamics with multiplicative covariance inﬂation},
  pdfauthor={K. Takeda and T. Sakajo}
}
\fi

\begin{document}

\maketitle

\begin{abstract}
Data assimilation is a method of uncertainty quantification to estimate the hidden true state by updating the prediction owing to model dynamics with observation data.
As a prediction model, we consider a class of nonlinear dynamical systems on Hilbert spaces including the two-dimensional Navier-Stokes equations and the Lorenz'63 and '96 equations.
For nonlinear model dynamics, the ensemble Kalman filter (EnKF) is often used to approximate the mean and covariance of the probability distribution with a set of particles called an ensemble.
In this paper, we consider a deterministic version of the EnKF known as the ensemble transform Kalman filter (ETKF), performing well even with limited ensemble sizes in comparision to other stochastic implementations of the EnKF.
When the ETKF is applied to large-scale systems, an ad-hoc numerical technique called a covariance inflation is often employed to reduce approximation errors.
Despite the practical effectiveness of the ETKF, little is theoretically known.
The present study aims to establish the theoretical analysis of the ETKF.
We obtain that the estimation error of the ETKF with and without the covariance inflation is bounded for any finite time.
In particular, the uniform-in-time error bound is obtained when an inflation parameter is chosen appropriately, justifying the effectiveness of the covariance inflation in the ETKF.
\end{abstract}



\section{Introduction}
\label{sec:introduction}
Data assimilation is a method of uncertainty quantification to estimate the hidden true state by combining model dynamics with observation data.
Data assimilation encompasses various problems such as interpolation of missing data by model dynamics.
In this paper, we consider the filtering problem \cite{lawDataAssimilationMathematical2015, reichProbabilisticForecastingBayesian2015}, where we estimate the time series of the true state using observations up to the current time.
We suppose that observation data is obtained by applying a linear operator to the true state and adding noise, which is the standard setting adopted in many filtering problems.
We then solve the problem with a Bayesian approach, in which the uncertainty in the state estimation is represented by a probability distribution.
One assimilation cycle consists of prediction and analysis steps.
In the prediction step, the dynamical model provides the prior distribution for the current state. 
It is then updated to the posterior distribution by Bayes' formula from the current observation in the analysis step.
The time series of the true state is estimated by repeating this assimilation cycle iteratively.
This approach is widely used in many problems.
See the comprehensive textbooks \cite{kalnayAtmosphericModelingData2002, reichProbabilisticForecastingBayesian2015}.

When the model dynamics is a linear system, the standard Bayesian assimilation method is the Kalman filter (KF) \cite{kalmanNewApproachLinear1960}, in which the probability distribution of the state is represented by the Gaussian distribution.
A remarkable property of linear dynamical systems is that they transform one Gaussian distribution to another.
Hence, we directly obtain the prior Gaussian distribution for the current time from that for the previous time in the prediction step.
Then, in the analysis step, the prior Gaussian distribution is updated to the posterior one using Bayes' formula with the current observation through the linear observation operator.
The KF is theoretically well-understood, and it works practically as well under the assumption that the additive observation noises and the probability distribution for the initial state follow Gaussian distributions \cite{lawDataAssimilationMathematical2015}.

In the meantime, for nonlinear model dynamics, the ensemble Kalman filter (EnKF) \cite{evensenDataAssimilationEnsemble2009} is employed as an extension of the KF.
Since the nonlinear dynamical system doesn't necessarily transfer one Gaussian distribution to another, we cannot assume the probability distribution is always represented by a Gaussian distribution.
Hence, the EnKF approximates probability distributions by a set of particles called an ensemble.
The well-known methods of the EnKF are the perturbed observation (PO) method \cite{burgersAnalysisSchemeEnsemble1998} and the ensemble transform Kalman filter (ETKF) \cite{bishopAdaptiveSamplingEnsemble2001}.
While both methods use the same prediction step in which each ensemble member evolves according to the nonlinear dynamical model, they take different approaches to update the ensemble in the analysis step.
The PO method adds Gaussian noises to the observation data in order to replicate artificial observations for each ensemble member.
Each member of the analysis ensemble is then obtained by the weighted average of a replicated observation and each member of the prediction ensemble.
When the ensemble size is large, this method performs effectively despite that artificial noises give rise to further uncertainties.
However, the PO method can be inaccurate with small ensemble sizes.
On the other hand, the ETKF generates the analysis ensemble deterministically.
The mean of the analysis ensemble is obtained by the weighted average of the observation and the prediction ensemble mean.
The deviations of the analysis ensemble are determined by the linear transformation of the prediction ensemble so that the covariance coincides with that obtained in the KF.
It is numerically confirmed that the ETKF performs well even with small ensemble sizes owing to its deterministic implementation \cite{bishopAdaptiveSamplingEnsemble2001}. 
Another practical advantage of the ETKF is that it is unnecessary to evaluate the ensemble covariance, avoiding redundant memory allocations in numerical computation \cite{huntEfficientDataAssimilation2007}.
Let us note that both methods are theoretically consistent with the KF \cite{kwiatkowskiConvergenceSquareRoot2015a,mandelConvergenceEnsembleKalman2011}.
That is to say, the ensemble mean and covariance converge to those in the KF in the limit of $ N \rightarrow \infty $ when a linear system with Gaussian noises is considered.

In the EnKF, the approximation of the covariance is often underestimated due to limited ensemble sizes, which causes an inaccurate estimation of the true state.
To handle the problem, ad-hoc covariance inflation techniques \cite{majdaFilteringComplexTurbulent2012a, tongNonlinearStabilityErgodicity2016a} are implemented in an additive or multiplicative way.
In an additive inflation,  one adds a positive diagonal matrix to the ensemble prediction covariance.
This is applied to the PO method, and it enhances the accuracy of the state estimation for a long time.
However, it is not applied to the ETKF since the covariance is not explicitly calculated in the method.
Instead, a multiplicative inflation is considered, in which one multiplies the covariance or ensemble deviation by a scaling factor.

While many nonlinear filtering problems are considered on finite-dimensional model dynamics \cite{dewiljesLongtimeStabilityAccuracy2018,lawDataAssimilationMathematical2015,tongNonlinearStabilityErgodicity2016a}, most dynamics such as the atmosphere and ocean are modeled by nonlinear and chaotic partial differential equations \cite{kalnayAtmosphericModelingData2002,reichProbabilisticForecastingBayesian2015}, which are infinite-dimensional dynamical systems.
Hence, it is theoretically essential to study nonlinear data assimilation methods in infinite dimensional spaces.
For the PO method, Kelly et al. \cite{kellyWellposednessAccuracyEnsemble2014b} shows that the filtering error does not grow faster than exponentially over time.
They also prove that an appropriate additive inflation ensures a uniform-in-time error bound, clarifying the effect of additive covariance inflation in the PO method.
On the other hand, the error bound of the ETKF remains unexplored due to the complexity of its elaborate algorithm.
The purpose of this paper is to establish the mathematical analysis of the ETKF for infinite-dimensional dynamical systems.
We reveal the basic properties of the ETKF with and without multiplicative covariance inflation.
This analysis theoretically validates the efficiency of the ETKF and the multiplicative covariance inflation.

This paper is constructed as follows.
In Section \ref{sec:setup}, we introduce some notations and the concept of nonlinear filtering problems, and explain the algorithm of the ETKF and the covariance inflation technique.
In Section \ref{sec:results}, we show the two main results for the error analysis of the ETKF.
Section \ref{sec:discussion} is a summary of the results and a discussion of future directions.

\section{Nonlinear data assimilation problem and algorithm}
\label{sec:setup}
\subsection{Notations}
Let $ \U $ be a separable Hilbert space endowed with the norm $ |\wc| $ and the inner product $ \bracket{\wc}{\wc} $, and $ \L(\U, \Y) $ denote the space of all bounded linear operators between two Hilbert spaces $ \U $ and $ \Y \subset \U $.
For $ N \in \N $, by $ I_\U $ and $ I_N $, we express the identity maps on $ \U $ and $ \R^N $ respectively.
For $ A \in \L(\U) \coloneq \L(\U, \U)$, $ \opnorm{A} $ represents the operator norm of $ A $, $ \ran(A) $ denotes the range of $ A $, and $ A^* $ is the adjoint of $ A $.
We also use $ u^* \in \U^* $ as a dual $ u \in \U $.
For $ u, v \in \U $, we define their product $ u \otimes v \in \L(\U) $ by $ u \otimes v: \U \ni w \mapsto u \bracket{v}{w} \in \U $, which is equivalent to $ u v^* = u \otimes v $.

Let $ \L_{sa}(\U) $ denote the space of all self-adjoint operators in $ \L(\U) $ (i.e., $ A^* = A $). 
For $ A \in \L_{sa}(\U) $, by $ A \succeq 0 $, we mean that $ A $ is positive semidefinite, i.e., $ \bracket{u}{Au} \ge 0 $ for all $ u \in \U $.
On the other hand, when there exists $ c > 0 $ such that $ \bracket{u}{Au} \ge c |u|^2 $ for all $ u \in \U $, we say $ A $ is positive define, denoted by $ A \succ 0 $.
For $ A, B \in \L_{sa}(\U) $, the order $ A \succ (\text{resp.} \succeq) \, B $ means $ A - B \succ (\text{resp.} \succeq) \, 0 $.

Let $ N \in \N $, for $ U = [u^{(n)}]_{n=1}^N $ and $ V = [v^{(n)}]_{n=1}^N \in \U^N $, the $\ell_2 $-norm $ \enorm{U} $ and the products $ U V^* \in \L(\U) $ and $ U^* V \in \R^{N \times N} $ are given by
\begin{align*}
    \enorm{U} = \left(\frac{1}{N} \sum_{n=1}^N \left|u^{(n)}\right|^2\right)^{\frac{1}{2}}, \quad 
    U V^* = \sum_{n=1}^N u^{(n)} \otimes v^{(n)}, \quad
    U^* V = \left[\bracket{u^{(i)}}{v^{(j)}}\right]_{i, j=1}^N.
\end{align*}
When we write $ u\bm{1} = [u, \dots, u] \in \U^N $ by $ \bm{1} = (1, \dots, 1) \in (\R^N)^* $, it holds that $ \enorm{u\bm{1}}^2 = |u|^2 $ for all $ u \in \U $.
Moreover, for $ m \in \U $ and $ T \in \R^{N \times N} $, we define
\begin{align*}
    m + U &= m \bm{1} + U = [m + u^{(n)}]_{n=1}^N \in \U^N, \quad UT = \left[\sum_{k=1}^N u^{(k)} T_{k,n}\right]_{n=1}^N \in \U^N.
\end{align*}
For an ensemble $ V = [v^{(n)}]_{n=1}^N \in \U^N $, $ \overline{v} = \frac{1}{N} \sum_{n=1}^N v^{(n)} $ is the ensemble mean and $ dV = [v^{(n)} - \overline{v}]_{n=1}^N \in \U^N $ is the ensemble deviation.
The ensemble $ V $ is then decomposed into the mean and the deviations, $ V = \overline{v}\bm{1} + dV $.
The (unbiased) ensemble covariance $ \cov(V) \in \L(\U) $ is defined by
\begin{align*}
    \cov(V) = \frac{1}{N-1} dV dV^*.
\end{align*}
Note that it is easy to see $ \cov(V) = \cov(dV) $ and $ \cov(V) \succeq 0 $.

The following lemma provides some equivalent representations of the $\ell_2 $-norm in terms of an ensemble $ V $.
\begin{lemma}
    \label{lem:l2_norm}
    The $ \ell_2 $-norm satisfies
    \begin{align}
        \label{eq:l2_norm}
        \enorm{V}^2 = \frac{1}{N} \trace{V^*V} = \frac{1}{N} \trace{VV^*} = |\overline{v}|^2 + \enorm{dV}^2.
    \end{align}
\end{lemma}
\begin{proof}
    The first equality is derived from the definition of $ \enorm{V} $.
    \begin{align*}
        \enorm{V}^2 = \frac{1}{N} \sum_{n=1}^N |v^{(n)}|^2 = \frac{1}{N} \sum_{n=1}^N \bracket{v^{(n)}}{v^{(n)}} = \frac{1}{N} \trace{V^*V}.
    \end{align*}
    Let $ (\phi_i)_{i \in \N} $ be a complete orthonormal basis of $ \U $, we have
    \begin{align*}
        \enorm{V}^2
        & = \frac{1}{N} \sum_{n=1}^N |v^{(n)}|^2 = \frac{1}{N} \sum_{n=1}^N \sum_{i \in \N} \bracket{v^{(n)}}{\phi_i}^2 = \frac{1}{N} \sum_{i \in \N} \sum_{n=1}^N \bracket{v^{(n)}}{\phi_i}^2 \\
        & = \frac{1}{N} \sum_{i \in \N} \sum_{n=1}^N \bracket{\phi_i}{(v^{(n)} \otimes v^{(n)}) \phi_i} = \frac{1}{N} \trace{VV^*}.
    \end{align*}
    Owing to the relation $ dV \bm{1}^* = 0 $, we have $ VV^* = \overline{v}\bm{1} \bm{1}^*\overline{v}^* + dVdV^* = N\overline{v}\,\overline{v}^* + dVdV^* $.
    Hence, we obtain $ \frac{1}{N} \trace{VV^*} = |\overline{v}|^2 + \enorm{dV}^2 $.
\end{proof}

\subsection{State space model}
We consider the following dynamical system on $\U $.
\begin{align}
    \label{eq:dynamics}
    \frac{du}{dt} = F(u), \quad u \in \U.
\end{align}
We suppose that a unique solution exists for any $ u_0 \in \U $ and it generates an analytic semigroup $ \Psi_t: \U \rightarrow \U $ for $ t \ge 0 $.
For a fixed period $ h > 0 $, we denote $ \Psi = \Psi_h $ and consider a discrete dynamical system 
\begin{align}
    \label{eq:dynamics_disc}
    u_j = \Psi(u_{j-1}), \quad j \in \N .
\end{align}
The Hilbert space $ \Y \subset \U $ is the observation space.
Then, the noisy observation $ y_j \in \Y $ is obtained from $ u_j $ for each $ j \in \N $,
\begin{align}
    \label{eq:obs}
    y_j = H u_j + \xi_j, \quad \xi_j \sim N(0, \Gamma),
\end{align}
where $ H \in \L(\U, \Y) $ is the observation operator  and $ \Gamma \in \L_{sa}(\Y) $ is the covariance operator with $ \Gamma \succ 0 $.
Remark that the assumption $ \Gamma \succ 0 $ is necessary to ensure the bounded inverse $ \Gamma^{-1} $, which makes the algorithm well-defined as discussed later.

Let $ (\Omega, \F, \Prob) $ be a probability space and $ \E $ denote expectation with respect to $ \Prob $.
We assume that the i.i.d.\,noise sequence $ (\xi_j)_{j \in \N} $ is independent of $ u_0 $.
The filtering problem aims to construct an approximation of $ (u_j)_{j \in \N} $ by a stochastic process $ (v_j)_{j \in \N} $ adapted to the filtration $ (\F_j)_{j \in \N} \subset \F $ where each $ \sigma $-algebra $ \F_j $ is generated by initial uncertainties $ u_0 $, $ v_0 $, and the noise sequence $ (\xi_k)_{k=1}^j $ for $ j \in \N $.
The conditional probability distribution $ \Prob^{v_j}(\wc|Y_j) $ is called by the filtering distribution for provided observations $ Y_j = (y_k)_{k=1}^j $.
Here, we focus on an algorithm computing $ \Prob^{v_j}(\wc|Y_j) $ from given $ \Prob^{v_{j-1}}(\wc|Y_{j-1}) $ and $ y_j $.
The algorithm consists of two steps, which we refer to as the prediction and analysis steps.
In the prediction step, the estimate of the current state $ \widehat{v}_j $ is provided by the time evolution of the dynamics $ \widehat{v}_j = \Psi(v_{j-1}) $.
Accordingly, the filtering distribution $ \Prob^{v_{j-1}}(\wc|Y_{j-1}) $ is propagated into the predicted distribution $ \Prob^{\widehat{v}_j}(\wc|Y_{j-1}) $.
In the analysis step, the current observation $ y_j $ is assimilated into the prediction $ \widehat{v}_j $ to obtain $ v_j $, and the predicted distribution $ \Prob^{\widehat{v}_j}(\wc|Y_{j-1}) $ is updated to $ \Prob^{v_j}(\wc|Y_j) $.

\subsection{The ensemble transform Kalman filter}
For $ N \in \N $, the ensemble transform Kalman filter (ETKF) approximates the filtering distribution by the empirical distribution of an ensemble of particles $ V_j = [v^{(1)}_j, \dots, v^{(N)}_j] \in \U^N $, $ \Prob^{v_j}(\wc|Y_j) \approx \frac{1}{N} \sum_{n=1}^N \delta_{v^{(n)}_j}(\wc) $, where $ \delta_a(\wc) $ denotes the Dirac measure. 
In the prediction step, each particle is driven by the dynamics \cref{eq:dynamics_disc}.
In the analysis step, the ensemble deviation is transformed by a matrix so that the posterior ensemble covariance corresponds to that in the minimum variance posterior ensemble.
A detailed description is given as follows.
\begin{definition}
    \label{def:etkf}
    The algorithm of the ETKF is as follows.
    \begin{enumerate}[label=(\Roman*)]
        \item {\label{def:etkf:1}Draw $ N $ independent samples $ V_0 = [v^{(n)}_0]_{n=1}^N $ from some probability distribution $ \mu_0 $.}
        \item {\label{def:etkf:2} (the prediction step, input: $ V_{j-1} \rightarrow $ output: $ \widehat{V}_j $) Compute the time evolution $ \widehat{v}_j^{(n)} = \Psi(v_{j-1}^{(n)}), \, n = 1, \dots, N $ and set $ \widehat{V}_j = [\widehat{v}_j^{(n)}]_{n=1}^N \in \U^N $.
        }
        \item {\label{def:etkf:3} (the analysis step, input: $ \widehat{V}_j, y_j \rightarrow $ output: $ V_j $) 
        Decompose $ \widehat{V}_j = \overline{\widehat{v}}_j\bm{1} + d\widehat{V}_j $, and put $ \widehat{C}_j = \cov(\widehat{V}_j) $.
        Update the mean
        \begin{align}
            \label{def:etkf:mean}
            \overline{v}_j = \overline{\widehat{v}}_j + K_j(y_j - H \overline{\widehat{v}}_j)
        \end{align} 
        with the Kalman gain,
        \begin{align}
            \label{def:kalman_gain}
            K_j = \widehat{C}_jH^* (H\widehat{C}_jH^* + \Gamma)^{-1}.
        \end{align}
        Take a symmetric transform matrix $ T_j \in \R^{N \times N} $ satisfying
        \begin{align}
            \label{def:etkf:transform}
            \frac{1}{N-1} d\widehat{V}_j T_j (d\widehat{V}_j T_j)^* = (I - K_j H) \widehat{C}_j,
        \end{align}
        and transform the ensemble deviation $ dV_j = d\widehat{V}_j T_j $.
        Finally, set the analysis ensemble $ V_j = \overline{v}_j\bm{1} + dV_j $.
        }
    \end{enumerate}
    \textnormal{Practically, we can use the following analysis step instead of \ref{def:etkf:3}.
    This avoids redundant memory allocations in practical numerical computation since we don't need to evaluate the covariance $ \widehat{C}_j $ explicitly.}
    \begin{itemize}
        \item[(III')\label{def:etkf:3'}] {(Input: $ \widehat{V}_j, y_j \rightarrow $ Output: $ V_j $)
        Decompose $ \widehat{V}_j = \overline{\widehat{v}}_j\bm{1} + d\widehat{V}_j $.
        Define the transform matrix
        \begin{align}
            \label{def:transform_matrix}
            T_j = \left(I_N + \frac{1}{N-1} d\widehat{V}_j^* H^* \Gamma^{-1} Hd\widehat{V}_j \right)^{-\frac{1}{2}},
        \end{align}
        and
        \begin{align*}
            \widetilde{T}_j = \frac{1}{N-1} T_j^2 d\widehat{V}_j^* H^* \Gamma^{-1} (y_j - H \overline{\widehat{v}}_j)\bm{1} + T_j.
        \end{align*}
        Transform
        \begin{align*}
            V_j = \overline{\widehat{v}}_j\bm{1} + d\widehat{V}_j \widetilde{T}_j.
        \end{align*}
        }
    \end{itemize}
\end{definition}
Note that \hyperref[def:etkf:3']{(III')} is equivalent to \ref{def:etkf:3}. 
The following theorem shows the existence of the transform matrix $ T_j $, which indicates that the ETKF is well-defined.
\begin{thm}
    \label{thm:welldefined}
    For any $ \widehat{V}_j \in \U^N $, there exists a unique symmetric transform matrix $ T_j \in \R^{N \times N} $ satisfying \cref{def:etkf:transform}.
\end{thm}
To prove \Cref{thm:welldefined}, let us summarize the key property of the Kalman gain $ K_j $.
\begin{lemma}
    \label{lem:kf}
    Let $ \widehat{C}_j \succeq 0 $, $ H \in \L(\U, \Y) $, and $ \Gamma \in \L(\Y) $ with $ \Gamma \succ 0 $.
    The Kalman gain satisfies
    \begin{align}
        \label{lem:kf:opti}
        K_j = (I_\U - K_jH)\widehat{C}_jH^*\Gamma^{-1},
    \end{align}
    and \cref{def:etkf:mean} is equivalent to 
    \begin{align}
        (I_\U + \widehat{C}_j H^* \Gamma^{-1} H)\overline{v}_j = \overline{\widehat{v}}_j + \widehat{C}_j H^* \Gamma^{-1} y_j. \label{eq:mean_update}
    \end{align}
\end{lemma}
\begin{proof}
    For simplicity, we omit the time index $ j $ in the following proofs since $ j \in \N $ is fixed.
    Owing to $ \widehat{C} \succeq 0 $ and $ \Gamma \succ 0 $, we have $ \Gamma + H\widehat{C}H^* \succ 0 $ and it is thus invertible.
    Then, $ I_\Y + \Gamma^{-1}H\widehat{C}H^* = \Gamma^{-1} (\Gamma + H\widehat{C}H^*) $ is also invertible since a product of two positive definite operators has positive spectrum \cite{hladnikSpectrumProductOperators1988}.
    From \cref{lem:i+a:inv} and the fact that $ (AB)^{-1} = B^{-1}A^{-1} $ for invertible $ A, B $, we have
    \begin{align*}
        (\Gamma + H\widehat{C}H^*)^{-1}
        & = (I_\Y + \Gamma^{-1}H\widehat{C}H^*)^{-1}\Gamma^{-1} = [I_\Y - (I_\Y + \Gamma^{-1}H\widehat{C}H^*)^{-1}\Gamma^{-1}H\widehat{C}H^*]\Gamma^{-1} \\
        & = [I_\Y - (\Gamma + H\widehat{C}H^*)^{-1}H\widehat{C}H^*]\Gamma^{-1}.
    \end{align*}
    Hence, we have.
    \begin{align*}
        K & = \widehat{C}H^*(\Gamma + H\widehat{C}H^*)^{-1} = \widehat{C}H^*[I_\Y - (\Gamma + H\widehat{C}H^*)^{-1}H\widehat{C}H^*]\Gamma^{-1} \\
        & = (I_\Y - KH) \widehat{C}H^* \Gamma^{-1},
    \end{align*}
    which is \cref{lem:kf:opti}.
    On the other hand, from \cref{lem:kf:opti}, we have 
    \begin{align*}
        \overline{v} = \overline{\widehat{v}} + K(y - H \overline{\widehat{v}}) = (I_\U - KH)\overline{\widehat{v}} + Ky =  (I_\U - KH)\overline{\widehat{v}} + (I_\U - KH) \widehat{C}H^* \Gamma^{-1}y.
    \end{align*}
    To show \cref{eq:mean_update}, it is sufficient to check $ (I_\U + \widehat{C}H^* \Gamma^{-1} H)(I_\U - KH) = I_\U$ with \cref{def:kalman_gain}.
    This is confirmed straightforwardly as follows.
    \begin{align*}
        & (I_\U + \widehat{C}H^* \Gamma^{-1} H)(I_\U - KH) \\
        & = I_\U + \widehat{C}H^* \Gamma^{-1} H - \widehat{C}H^*(\Gamma + H\widehat{C}H^*)^{-1}H - \widehat{C}H^* \Gamma^{-1} H \widehat{C}H^*(\Gamma + H\widehat{C}H^*)^{-1}H \\
        & = I_\U + \widehat{C}H^* \Gamma^{-1} H - \widehat{C}H^*\Gamma^{-1}\Gamma(\Gamma + H\widehat{C}H^*)^{-1}H - \widehat{C}H^* \Gamma^{-1} H \widehat{C}H^*(\Gamma + H\widehat{C}H^*)^{-1}H \\
        & = I_\U + \widehat{C}H^* \Gamma^{-1} H - \widehat{C}H^*\Gamma^{-1} (\Gamma + H\widehat{C}H^*)(\Gamma + H\widehat{C}H^*)^{-1} H = I_\U.
    \end{align*}
\end{proof}

\begin{proof}[Proof of \Cref{thm:welldefined}]
    First, we prove the existence of $ T_j $ satisfying \cref{def:etkf:transform}.
    Let $ dY = H d\widehat{V} $.
    Then, the operator $ \Gamma + \frac{1}{N-1}dYdY^* = \Gamma + H\widehat{C}H^* $ is invertible, and we consider the symmetric matrix
    \begin{align*}
        S = I_N - \frac{1}{N-1} dY^* \left(\Gamma + H\widehat{C}H^*\right)^{-1}dY \in \R^{N \times N}.
    \end{align*}
    Then, we have
    \begin{align*}
        \frac{1}{N-1} d\widehat{V} S d\widehat{V}^*
        & = \frac{1}{N-1} d\widehat{V} d\widehat{V}^* - \frac{1}{N-1}d\widehat{V} dY^* \left(\Gamma + H\widehat{C}H^*\right)^{-1}\frac{1}{N-1} dY d\widehat{V}^* \notag \\
        & = \widehat{C} - \widehat{C}H^*(\Gamma + H\widehat{C}H^*)^{-1} H \widehat{C} = (I - KH)\widehat{C}, 
    \end{align*}
    in which, we use $ \frac{1}{N-1} dY d\widehat{V}^* = H \widehat{C} $.
    From \cref{lem:invertible_operator:main}, we have 
    \begin{align}
        \label{eq:s}
        S = (I_N + \frac{1}{N-1} dY^* \Gamma^{-1} dY)^{-1}
    \end{align}
    and $ S \succ 0 $.
    We finally define the transform matrix $ T = S^{\frac{1}{2}} $, which is nothing but \cref{def:transform_matrix}.
    Then $ T $ becomes symmetric by definition.
\end{proof}

We finally show another representation of the ensemble update.
\begin{lemma}
    \label{lem:etkf}
    The following hold for the ensemble deviation and the transform matrix
    \begin{align}
        dV_j \bm{1}^* & = d\widehat{V}_j \bm{1}^*  = 0 \in \U, \label{lem:etkf:ens_dev} \\
        T_j \bm{1}^* & = \bm{1}^*.\label{lem:etkf:id_on_1}
    \end{align}
    Moreover, the ensembles satisfy the relation
    \begin{align}
        (I_\U + \widehat{C}_j H^*\Gamma^{-1}H) V_j = \widehat{V}_j T_j^{-1} + \widehat{C}_j H^* \Gamma^{-1} y_j \bm{1} \label{eq:transform_repr}.
    \end{align}
\end{lemma}
\begin{proof}
    The definition of the ensemble mean yields the first property \cref{lem:etkf:ens_dev}.
    \begin{align*}
        dV \bm{1}^* = \sum_{n=1}^N (v^{(n)} - \overline{v}) = 0,
    \end{align*}
    yielding
    \begin{align*}
        S^{-1} \bm{1}^* = \left(I_N + \frac{1}{N-1} dV^* H^* \Gamma^{-1} H dV\right) \bm{1}^* = \bm{1}^*.
    \end{align*}
    Hence, we have
    \begin{align}
        \label{lem:etkf:eqS}
        S\bm{1}^* = \bm{1}^*.
    \end{align}
    Then, we prove that $ \bm{1}^* $ is also an eigenvector of $ T = S^{\frac{1}{2}} $ with an eigenvalue $ 1 $.
    Since $ S $ is symmetric, it is diagonalized as $ S = U D U^* $ for a unitary $ U \in \R^{N \times N} $ and a diagonal $ D \in \R^{N \times N} $.
    Then, \cref{lem:etkf:eqS} is equivalent to
    \begin{align*}
        S\bm{1}^* = \bm{1}^* \Leftrightarrow U D U^* \bm{1}^* = \bm{1}^* \Leftrightarrow D U^* \bm{1}^* = U^* \bm{1}^*.
    \end{align*}
    Putting $ \bm{u} = U^* \bm{1}^* = (u_1, \cdots, u_N)^* \in \R^N $ and $ d_n > 0 $ as $ n $-diagonal element of $ D $ for $ n = 1, \dots, N $, the last equality is rewritten for each component 
    \begin{align*}
        d_n u_n = u_n, \quad n = 1, \dots, N.
    \end{align*}
    This implies that $ d_n = 1 $ or $ u_n = 0 $ for each $ n = 1, \dots, N $.
    Hence, the following also holds
    \begin{align*}
        d_n^{\frac{1}{2}} u_n = u_n, \quad n = 1, \dots, N,
    \end{align*}
    and we have $ D^{\frac{1}{2}} U^* \bm{1}^* = U^* \bm{1}^* $.
    By definition, $ T $ is written as $ T = U D^{\frac{1}{2}} U^* $, and this yields $ T \bm{1}^* = \bm{1}^* $, which is \cref{lem:etkf:id_on_1}.

    The last equality \cref{eq:transform_repr} is shown as follows.
    From \eqref{eq:mean_update}, we have
    \begin{align*}
        (I_\U + \widehat{C}H^* \Gamma^{-1} H)\overline{v} = \overline{\widehat{v}} + \widehat{C}H^* \Gamma^{-1} y \in \U.
    \end{align*}
    By using $ \widehat{C} = \cov(\widehat{V}) $, \cref{eq:s} and $ S = T^2 $, we obtain
    \begin{align*}
        (I_\U + \widehat{C}H^* \Gamma^{-1} H) dV
        & = (I_\U + \widehat{C}H^* \Gamma^{-1} H) d\widehat{V} T = d\widehat{V} \left[I_\U + \frac{1}{N-1} d\widehat{V}^* H^* \Gamma^{-1} H d\widehat{V}\right]T \\
        & = d\widehat{V} S^{-1} T = d\widehat{V} T^{-1} \in \U^N.
    \end{align*}
    Finally, owing to \cref{eq:mean_update} and $ \overline{\widehat{v}}\bm{1} T^{-1} = \overline{\widehat{v}}\bm{1} $,
    \begin{align*}
        (I_\U + \widehat{C}H^* \Gamma^{-1} H) V
        & = (I_\U + \widehat{C}H^* \Gamma^{-1} H) (\overline{v}\bm{1} + dV) = \overline{\widehat{v}}\bm{1} + \widehat{C}H^* \Gamma^{-1} y\bm{1} + d\widehat{V} T^{-1} \\
        & = \overline{\widehat{v}}\bm{1} T^{-1} + d\widehat{V} T^{-1} + \widehat{C}H^* \Gamma^{-1} y\bm{1} = \widehat{V} T^{-1} + \widehat{C}H^* \Gamma^{-1} y\bm{1}.
    \end{align*}
    This finishes the proof.
\end{proof}

\section{Bounds for the filtering error of the ensemble transform Kalman filter}
\label{sec:results}
For each $ j \in \N \cup \{0\} $, we define the filtering ensemble error by
\begin{align*}
    E_j = [e_j^{(n)}]_{n=1}^N \quad \text{with} \quad e_j^{(n)} & = v_j^{(n)} - u_j, \quad n = 1, \dots N.
\end{align*}
Note that $ dE_j = dV_j $, and hence, $ E_j = e_j\bm{1} + dV_j $ where $ e_j = \overline{v}_j - u_j \in \U $.
Furthermore, owing to \cref{eq:l2_norm}, the decomposition of the norm of the ensemble error is given by
\begin{align}
    \label{eq:error_decomp}
    \enorm{E_j}^2 = |e_j|^2 + \enorm{dV}^2.
\end{align}

\subsection{Well-posedness}
\label{subsec:wellposed}
We show the well-posedness of the ETKF, i.e., the ensemble error does not blow up faster than an exponential function.
To do this, we make some assumptions on the dynamics and the observation to analyze the properties of the ETKF as in \cite{kellyWellposednessAccuracyEnsemble2014b}.
\begin{assumption}
    \label{asmp:model_base}
    There exists $ \rho > 0 $ such that $ \Psi_t $ has an absorbing ball $ B(\rho) = \{v \in \U \mid |v| \le \rho \} $, i.e., $ \Psi_t(v) \in B(\rho) $ for any $ v \in B(\rho) $ and $ t \ge 0 $.
\end{assumption}

\begin{assumption}
    \label{asmp:model1}
    There exists $ \beta > 0 $ such that, for any $ u \in B(\rho) $ and $ v \in \U $,
    \begin{align}
        \bracket{F(u) - F(v)}{u-v} \le \beta |u-v|^2.
    \end{align}
\end{assumption}
Note that \Cref{asmp:model_base,asmp:model1} are satisfied by Lorenz '63, Lorenz '96 models and the incompressible Navier-Stokes equations on a two-dimensional torus \cite{majdaNonlinearDynamicsStatistical2006,majdaFilteringComplexTurbulent2012a,temamInfinitedimensionalDynamicalSystems1997}.
The following lemma was established in \cite{kellyWellposednessAccuracyEnsemble2014b}, providing the upper bound of the error growth rate due to the model dynamics.
\begin{lemma}
    \label{lem:chaos}
    Suppose that \cref{asmp:model_base,asmp:model1} hold.
    Then, for any $ u \in B(\rho) $, $ v \in \U $ and $ t > 0 $,
    \begin{align}
        |\Psi_t(u) - \Psi_t(v)| \le e^{\beta t}|u-v|.
    \end{align}
\end{lemma}
\begin{proof}
    See Lemma 2.6 of \cite{kellyWellposednessAccuracyEnsemble2014b}.
\end{proof}

Regarding the observation operator and its covariance operator, we assume that they are trivial.
\begin{assumption}
    \label{asmp:obs1}
    The state is fully observed, i.e., $ H =  I_\U $ and $ \Gamma = \gamma^2 I_\U $ for some $ \gamma > 0 $.
\end{assumption}

\begin{thm}
    \label{thm:wellposed}
    Under Assumptions \ref{asmp:model_base}--\ref{asmp:obs1}, we consider the ETKF provided in \Cref{def:etkf}.
    Then, we have the following upper bound.
    \begin{align}
        \label{thm:wellposed:ineq}
        \E\left[\enorm{E_j}^2\right] \le e^{2\beta h j} \E\left[\enorm{E_0}^2\right] + (N-1)\gamma^2 \frac{e^{2\beta h j} - 1}{e^{2\beta h} - 1}, \quad j \in \N.
    \end{align}
\end{thm}
\begin{proof}
    From \Cref{asmp:obs1}, the relation \eqref{eq:transform_repr} becomes
    \begin{align}
        \label{thm:wellposed:pf:eq0}
        (I_\U + \gamma^{-2}\widehat{C}_j) V_j = \widehat{V}_jT^{-1}_j + \gamma^{-2}\widehat{C}_j y_j \bm{1}.
    \end{align}
Let $ U_j = u_j\bm{1} \in \U^N $.
From \cref{lem:etkf:id_on_1}, we have $ U_j = U_j T_j^{-1} $.
Hence,
\begin{align}
    \label{thm:wellposed:pf:eq1}
    (I_\U + \gamma^{-2}\widehat{C}_j) U_j = U_j + \gamma^{-2}\widehat{C}_j U_j = U_j T_j^{-1} + \gamma^{-2}\widehat{C}_j U_j.
\end{align}
Setting $ \widehat{E}_j = \widehat{V}_j - U_j $ and subtracting \cref{thm:wellposed:pf:eq1} from \cref{thm:wellposed:pf:eq0} yields
\begin{align*}
    (I_\U + \gamma^{-2}\widehat{C}_j) E_j = \widehat{E}_j T_j^{-1} + \gamma^{-2}\widehat{C}_j (y_j - u_j)\bm{1} = \widehat{E}_j T_j^{-1} + \gamma^{-2}\widehat{C}_j \xi_j \bm{1}.
\end{align*}
Owing to $ \gamma^{-2}\widehat{C}_j \succeq 0 $, $ I_\U + \gamma^{-2}\widehat{C}_j $ is invertible.
Multiplying $ (I_\U + \gamma^{-2}\widehat{C}_j)^{-1} $, we obtain
\begin{align*}
    E_j = (I_\U + \gamma^{-2}\widehat{C}_j)^{-1} \widehat{E}_j T_j^{-1} + (I_\U + \gamma^{-2}\widehat{C}_j)^{-1}\gamma^{-2}\widehat{C}_j \xi_j \bm{1}.
\end{align*}
Let us divide $ E_j $ into the following two terms and evaluate them separately.
\begin{align}
    R_1 & = (I_\U + \gamma^{-2}\widehat{C}_j)^{-1} \widehat{E}_j T_j^{-1}, \label{def:r1}\\
    R_2 & = (I_\U + \gamma^{-2}\widehat{C}_j)^{-1} \gamma^{-2}\widehat{C}_j \xi_j \bm{1} \label{def:r2}.
\end{align}

Here, the dimension of $ \ran(C_j) $ is $ N-1 $  at most since $ \widehat{C}_j $ consists of $ N $ vectors with one constraint. 
Let $ \Pi_j $ be the projection to $ \ran(C_j) $, and we have
\begin{align*}
    R_2 & = (I_\U + \gamma^{-2}\widehat{C}_j)^{-1} \gamma^{-2}\widehat{C}_j \Pi_j \xi_j \bm{1}.
\end{align*}
From \cref{lem:i+a:succeq}, we have $ (I_\U + \gamma^{-2}\widehat{C}_j)^{-1} \gamma^{-2}\widehat{C}_j \preceq I $ owing to $ \gamma^{-2}\widehat{C}_j \succeq 0 $.
This leads to
\begin{align}
    \label{ineq:r2_estimate}
    \enorm{R_2}^2 \le \enorm{\Pi_j\xi_j \bm{1}}^2 = |\Pi_j \xi_j|^2.
\end{align}
Let $ J = I_\U + \gamma^{-2}\widehat{C}_j $.
Then we have $ J, J^{-1} \in \L_{sa}(\U) $ and $|J^{-1}| \le 1 $.
We obtain
\begin{align*}
    R_1 R_1^* = J^{-1} \widehat{E}_j T_j^{-2} \widehat{E}_j^*  J^{-1}.
\end{align*}
Considering the relations $ d\widehat{E}_j = d\widehat{V}_j $ and $ \widehat{E}_jd\widehat{V}_j^* = d\widehat{V}_j\widehat{E}_j^* = d\widehat{V}_j d\widehat{V}_j^* $,
we have
\begin{align*}
    \widehat{E}_j T_j^{-2} \widehat{E}_j^*
    & = \widehat{E}_j \left[I_N + \frac{\gamma^{-2}}{N-1} d\widehat{V}_j^* d\widehat{V}_j\right] \widehat{E}_j^* \\
    & = \widehat{e}_j\bm{1} (\widehat{e}_j\bm{1})^* + d\widehat{V}_jd\widehat{V}_j^* + d\widehat{V}_jd\widehat{V}_j^* \gamma^{-2} \widehat{C}_j\\
    & = \widehat{e}_j\bm{1} (\widehat{e}_j\bm{1})^* + d\widehat{V}_jd\widehat{V}_j^* J,
\end{align*}
where $ \widehat{e}_j = \overline{\widehat{v}}_j - u_j $.
Since $ J^{-1} $ is self-adjoint, we have
\begin{align*}
    R_1 R_1^* = J^{-1} \widehat{e}_j\bm{1} (\widehat{e}_j\bm{1})^* J^{-1} + J^{-1} d\widehat{V}_jd\widehat{V}_j^* = J^{-1} \widehat{e}_j\bm{1} (J^{-1}\widehat{e}_j\bm{1})^* + J^{-1} d\widehat{V}_jd\widehat{V}_j^*.
\end{align*}
Then, $ \enorm{R_1}^2 $ is bounded by
\begin{align*}
    \enorm{R_1}^2 
    & = |J^{-1}\widehat{e}_j|^2 + \frac{1}{N}\trace{J^{-1} d\widehat{V}_jd\widehat{V}_j^*} \le |J^{-1}|^2|\widehat{e}_j|^2 + |J^{-1}| \frac{1}{N}\trace{d\widehat{V}_jd\widehat{V}_j^*} \\
    & \le |\widehat{e}_j|^2 + \enorm{d\widehat{V}_j}^2 = \enorm{\widehat{E}_j}^2.
\end{align*}
The first inequality follows from the H\"{o}lder inequality about the trace norm and the operator norm \cite{conwayCourseFunctionalAnalysis2007}, and the second inequality is owing to $ |J^{-1}| \le 1 $.
From this with \Cref{lem:chaos}, we obtain the upper bound of $ \enorm{R_1} $ as follows.
\begin{align}
    \label{ineq:r1_estimate}
    \enorm{R_1}^2 \le \enorm{\widehat{E}_j}^2 
    \le e^{2\beta h}\enorm{E_{j-1}}^2.
\end{align}

Since $ R_1 $ and $ R_2 $ is conditionally independent under $ \F_{j-1} $, it follows from \eqref{ineq:r2_estimate} and \eqref{ineq:r1_estimate} that
\begin{align*}
    \E_{j-1}[\enorm{E_j}^2]
    & = \E_{j-1}[\enorm{R_1}^2] + \E_{j-1}[\enorm{R_2}^2] \le e^{2\beta h}\E_{j-1}[\enorm{E_{j-1}}^2] + \E_{j-1}[|\Pi_j \xi_j|^2] \\
    & = e^{2\beta h}\E_{j-1}[\enorm{E_{j-1}}^2] + (N-1)\gamma^2,
\end{align*}
where the conditional expectation is denoted by $ \E_{j-1}[\wc] \coloneq \E[\,\wc \mid \F_{j-1}] $.
Here, the conditional expectation satisfies $ \E[\E_{j-1}[\wc]] = \E[\wc] $ in general.
Therefore, taking the expectation yields
\begin{align*}
    \E[\enorm{E_j}^2] \le e^{2\beta h}\E[\enorm{E_{j-1}}^2] + (N-1)\gamma^2.
\end{align*}
Applying this inequality repeatedly, we obtain \eqref{thm:wellposed:ineq}.
\end{proof}

\subsection{Error analysis with multiplicative covariance inflation}
\label{subsec:error_analysis}
If the ensemble size $ N $ is smaller than the dimension of $ \U $, the ensemble covariance $ \widehat{C}_j $ degenerates.
This implies that $ \opnorm{J^{-1}} = \opnorm{(I_\U + \gamma^{-2}\widehat{C}_j)^{-1}} = 1 $.
In the PO method \cite{kellyWellposednessAccuracyEnsemble2014b}, a diagonal matrix to $ \widehat{C}_j $ is added to avoid this degeneration.
In the meantime, this additive inflation technique is not applicable to the ETKF because $ \widehat{C}_j $ is not explicitly used in the computation of the transform matrix \cref{def:transform_matrix}.
Hence, multiplicative covariance inflation is used, in which we multiply an inflation factor $ \alpha \ge 1 $ by the deviations of the ensemble $ d\widehat{V}_j $ after the prediction step \ref{def:etkf:2} of \Cref{def:etkf}.
To be specific, for the prediction ensemble $ \widehat{V}_j = \overline{\widehat{v}}_j + d\widehat{V}_j $, we define the inflated prediction ensemble $ \widehat{V}_j^\alpha $ by
\begin{align*}
    \widehat{V}_j^\alpha = \overline{\widehat{v}}_j + \alpha d\widehat{V}_j.
\end{align*}
In the analysis step \ref{def:etkf:3} of \Cref{def:etkf}, we use $ \widehat{V}_j^\alpha $ instead of $ \widehat{V}_j $.
First, it follows that $ \widehat{C}_j^\alpha = \cov(\widehat{V}_j^\alpha) = \alpha^2 \widehat{C}_j $.
Then, the mean update \cref{def:etkf:mean} is given by
\begin{align*}
    \overline{v}_j^\alpha = \overline{\widehat{v}}_j + K_j^\alpha(y_j - H \overline{\widehat{v}}_j),
\end{align*}
where the Kalman gain is $ K_j^\alpha = \widehat{C}_j^\alpha H^* (H\widehat{C}_j^\alpha H^* + \Gamma)^{-1} $.
Second, the transform matrix is given by
\begin{align*}
    T_j^\alpha = \left(I_N + (d\widehat{V}_j^\alpha)^* H^* \Gamma^{-1} H d\widehat{V}_j^\alpha \right)^{-\frac{1}{2}} = \left(I_N + \alpha^2 d\widehat{V}_j^* H^* \Gamma^{-1} H d\widehat{V}_j\right)^{-\frac{1}{2}},
\end{align*}
and the deviations are transformed by
\begin{align*}
    dV_j^\alpha = d\widehat{V}_j^\alpha T_j^\alpha,
\end{align*}
and we define $ V^\alpha = \overline{v}_j^\alpha + dV_j^\alpha $.

The relation between the prediction and analysis ensemble is summarized as follows.
The mean update \eqref{eq:mean_update} is given by
\begin{align}
    \label{eq:mean_update_inflation}
    (I_\U + \alpha^2 \widehat{C}_j H^* \Gamma^{-1} H)\overline{v}_j^\alpha = \overline{\widehat{v}}_j + \alpha^2 \widehat{C}_jH^* \Gamma^{-1} y_j,
\end{align}
and the analysis covariance satisfies
\begin{align}
    \label{eq:cov_update_inflation}
    C_j^\alpha = \frac{\alpha^2}{N-1} d\widehat{V}_j (I_N + \alpha^2 d\widehat{V}_j^* H^* \Gamma^{-1} H d\widehat{V}_j)^{-1} d\widehat{V}_j^*.
\end{align}
In what follows, we omit the superscript $ \alpha $ in the notations, since the changes in terms of $ \alpha $ explicitly appear as the multiplicative factor $ \alpha^2 $ in \cref{eq:mean_update_inflation,eq:cov_update_inflation}, and no confusion occurs.

It is known that the multiplicative inflation is insufficient to resolve the issue of the covariance degeneration when the state space is infinite-dimensional since it does not affect the eigenvectors of the covariance.
Hence, we need the following assumption that the actual dimension of the dynamics becomes finite and $ F $ is Lipschitz, which is stronger than \Cref{asmp:model1}.

\begin{assumption}
    \label{asmp:model_lip}
    The dimension of the state space $ \U $ is finite $ m \in \N $.
    For $ \rho > 0 $ in \Cref{asmp:model_base} and $ F $ in \cref{eq:dynamics}, there exists $ \beta > 0 $ such that for any $ u, v \in B(\rho) $
    \begin{align*}
        |F(u) - F(v)| \le \beta |u - v|.
    \end{align*}
\end{assumption}
Let $ u_0 \in B(\rho) $ and $ v_0^{(n)} \in B(\rho) $ for $ n = 1, \dots, N $.
We write $ \overline{v}_t = \frac{1}{N} \sum_{n=1}^N \Psi_t(v^{(n)}_0) $ and $ u_t = \Psi_t(u_0) $.
Then, we have the following lemma.
\begin{lemma}
    \label{lem:chaos_proj}
    Suppose \Cref{asmp:model_base,asmp:model_lip}.
    Then, for any $ \epsilon > 0 $, $ t > 0 $,
    \begin{align}
        |\overline{v}_t - u_t|^2 \le e^{2(\beta + \epsilon) t}(|\overline{v}_0 - u_0|^2 + D) - D,
    \end{align}
    where $ D = \frac{2\beta^2\rho^2}{2(\beta+\epsilon)\epsilon} $.
\end{lemma}
\begin{proof}
    We write $ e_t = \overline{v}_t - u_t $ and $ \overline{F}_t = \frac{1}{N} \sum_{n=1}^N F(v^{(n)}_t) $.
    Then, we have $ \frac{d}{dt}\overline{v}_t = \overline{F}_t $.
    Hence, we obtain
    \begin{align}
        \frac{1}{2}\frac{d}{dt} | e_t|^2
        & = \bracket{\frac{d}{dt} e_t}{ e_t} = \bracket{\overline{F}_t - F(u_t)}{e_t} \notag \\
        & \le |\overline{F}_t - F(u_t)|| e_t| \notag \\
        & \le \left(|F(\overline{v}_t) - F(u_t)| + |\overline{F}_t - F(\overline{v}_t)|\right)| e_t| \notag \\
        & \le \left(\beta | e_t| + |\overline{F}_t - F(\overline{v}_t)|\right)| e_t|, \label{lem:chaos_proj:ineq1}
    \end{align}
    where we can use \Cref{asmp:model_lip} owing to $ \overline{v}_t, u_t \in B(\rho) $.
    The second term is estimated by 
    \begin{align*}
        |\overline{F}_t - F(\overline{v}_t)|^2 \le \frac{N}{N^2} \sum_{n=1}^N |F(v^{(n)}_t) - F(\overline{v}_t)|^2 \le \frac{1}{N} \sum_{n=1}^N \beta^2 |v^{(n)}_t - \overline{v}_t|^2 \le 4 \beta^2 \rho^2.
    \end{align*}
    Substituting it into \cref{lem:chaos_proj:ineq1} and using the Young inequality, we have
    \begin{align*}
        \frac{d}{dt} | e_t|^2 
        & \le 2\beta | e_t|^2 + 4\beta\rho| e_t| = 2\beta | e_t|^2 + (2\beta\rho\epsilon^{-1/2}) (2\epsilon^{1/2}| e_t|) \\
        & \le 2\beta | e_t|^2 + 2\beta^2\rho^2 \epsilon^{-1} + 2\epsilon| e_t|^2 = 2(\beta + \epsilon)(| e_t|^2 + D),
    \end{align*}
    for any $ \epsilon > 0 $.
    Therefore, we obtain $ | e_t|^2 \le e^{2(\beta + \epsilon) t}(| e_0|^2 + D) - D $ from the Gronwall inequality.
\end{proof}

We now consider the error between the ensemble mean and the true state.
\begin{align*}
    e_j = \overline{v}_j - u_j \in \U.
\end{align*}
Here, $ \lambda_{min}(C) $ denotes the minimum eigenvalue of a covariance operator $ C $.
\begin{thm}
    \label{thm:error_bound}
    Under \Cref{asmp:model_base,asmp:model_lip,asmp:obs1}, we consider the ETKF with the multiplicative inflation with $ \alpha \ge 1 $.
    Suppose also that the ensemble size $ N \in \N $ is large enough to satisfy $ \lambda_{min}( C_0 ) \ge \lambda_0 $ with $ \lambda_0 > 0 $, and that $ v^{(n)}_j \in B(\rho) $ for all $ n = 1, \dots, N $ and $ j \in \N $.
    Then, for any $ \epsilon > 0 $, there exists $ \alpha_0 = \alpha_0(\rho, \beta, N, \lambda_0, \gamma, \epsilon) \ge 1 $ such that the following hold for any $ \alpha \ge \alpha_0 $.
    \begin{enumerate}[label=(\roman*)]
        \item {There exists $ \lambda_* = \lambda_*(\rho, \beta, N, \lambda_0, \gamma, \alpha) > 0 $ such that $ \lambda_{min}( \widehat{C}_j ) > \lambda_* $ for all $ j \in \N $.}
        \item {For $ j \in \N $ and $ \theta = (1+\frac{\alpha^2}{\gamma^2} \lambda_*)^{-2} e^{2(\beta + \epsilon) h} $,
            \begin{align}
                \label{ineq:main:result}
                \E[|e_j|^2] \le \theta^j (\E[|e_0|^2] + D) + m \gamma^2 \frac{1 - \theta^j}{1-\theta} + \left(\frac{(1 - \theta^j)(1 - \Theta)}{1 - \theta} - 1\right)D,
            \end{align}
            where $ D = \frac{2\beta^2\rho^2}{2(\beta+\epsilon)\epsilon} $ and $ \Theta = (1+\frac{\alpha^2_0}{\gamma^2} \lambda_*)^{-2}$.
            Moreover, if $ \theta < 1 $, we have
            \begin{align}
                \label{ineq:main:result_limsup}
                \limsup_{j \rightarrow \infty} \E[|e_j|^2] \le \frac{m \gamma^2}{1-\theta} + \left(\frac{1 - \Theta}{1 - \theta} - 1\right)D.
            \end{align}
        }
    \end{enumerate}
\end{thm}

\begin{proof}
    For convenience, we write $ \widehat{\lambda}^{min}_j = \lambda_{min}( \widehat{C}_j ) $ and $ \lambda^{min}_j = \lambda_{min}( C_j ) $.
    We first estimate the change from $ \lambda^{min}_{j-1} $ to $ \widehat{\lambda}^{min}_j $ in the prediction step.
    To this end, we write $ \widehat{v}^{(n)}_t = \Psi_t(v^{(n)}_{j-1}) $, $ \widehat{C}_t = \frac{1}{N-1} \sum_{n=1}^N (\widehat{v}^{(n)}_t - \overline{\widehat{v}}_t) \otimes (\widehat{v}^{(n)}_t - \overline{\widehat{v}}_t) $, $ \lambda_t = \lambda_{min}( \widehat{C}_t ) $ for $ t \in [0, h] $.
    The differentiation of $ \widehat{C}_t $ with respect to $ t $ yields
    \begin{align}
        \label{eq:derivative_cov}
        \frac{d}{dt}\widehat{C}_t = \frac{1}{N-1} \sum_{n=1}^N (F(\widehat{v}^{(n)}_t) - \overline{F}_t) \otimes (\widehat{v}^{(n)}_t - \overline{\widehat{v}}_t) + (\widehat{v}^{(n)}_t - \overline{\widehat{v}}_t) \otimes (F(\widehat{v}^{(n)}_t) - \overline{F}_t),
    \end{align}
    where $ \overline{F}_t = \frac{1}{N} \sum_{n=1}^N F(\widehat{v}^{(n)}_t) $.
    Owing to $ \frac{1}{N} \sum_{n=1}^{N} \widehat{v}^{(n)}_t - \overline{\widehat{v}}_t = 0 $, we have
    \begin{align}
        \label{eq:derivative_lambda}
        \frac{d}{dt} \widehat{C}_t = \frac{1}{N-1} \sum_{n=1}^N (F(\widehat{v}^{(n)}_t) - F(\overline{\widehat{v}}_t)) \otimes (\widehat{v}^{(n)}_t - \overline{\widehat{v}}_t) + (\widehat{v}^{(n)}_t - \overline{\widehat{v}}_t) \otimes (F(\widehat{v}^{(n)}_t) - F(\overline{\widehat{v}}_t)).
    \end{align}
    With a similar argument on the derivative of the eigenvalue of the covariance matrix in Section~3 of \cite{dewiljesLongtimeStabilityAccuracy2018}, there exists $ w \in \U $ with $ |w| = 1 $ such that
    \begin{align*}
        \frac{d}{dt} \lambda_t = \bracket{w}{ \frac{d}{dt}\widehat{C}_t  w}.
    \end{align*}
    To derive the lower bound $ \frac{d}{dt} \lambda_t $, we consider the absolute value of the right-hand side of \cref{eq:derivative_lambda}.
    Owing to $ |w| = 1 $ and \Cref{asmp:model_lip}, we have
    \begin{align*}
        \left|\bracket{w}{ \frac{d}{dt}\widehat{C}_t  w}\right|
        & \le \left|\frac{2}{N-1}\sum_{n=1}^N \bracket{F(\widehat{v}^{(n)}_t) - F(\overline{\widehat{v}}_t)}{w}\bracket{\widehat{v}^{(n)}_t - \overline{\widehat{v}}_t}{w}\right| \\
        & \le 2\left(\frac{1}{N-1}\sum_{n=1}^N \bracket{F(\widehat{v}^{(n)}_t) - F(\overline{\widehat{v}}_t)}{w}^2\right)^{\frac{1}{2}}\left(\frac{1}{N-1}\sum_{n=1}^N \bracket{\widehat{v}^{(n)}_t - \overline{\widehat{v}}_t}{w}^2\right)^{\frac{1}{2}} \\
        & \le \beta \frac{1}{N-1} \sum_{n=1}^{N} |\widehat{v}^{(n)}_t - \overline{\widehat{v}}_t|^2 \le 8\frac{N}{N-1} \beta \rho^2.
    \end{align*}
    The last inequality comes from owing to $ \widehat{v}^{(n)}_t \in B(\rho) $ by \Cref{asmp:model_base} and the assumption of \Cref{thm:error_bound}.
    Hence, for $ a = 8\frac{N}{N-1} \beta \rho^2 > 0 $,
    \begin{align*}
        \frac{d}{dt}\lambda_t \ge - a.
    \end{align*}
    Integrating it from $ 0 $ to $ t = h $, we have
    \begin{align}
        \label{ineq:lambda_min:chaos}
        \widehat{\lambda}^{min}_j = \lambda_h \ge e^{-ah} \lambda_0 = e^{-ah} \lambda^{min}_{j-1}.
    \end{align}

    The next step is to address the change in the eigenvalue in the analysis step. 
    From \Cref{asmp:obs1}, we have
    \begin{align*}
        C_{j-1} = \frac{\alpha^2}{N-1} d\widehat{V}_{j-1} (I_N + \alpha^2 \gamma^{-2}\widetilde{C}_{j-1})^{-1} d\widehat{V}_{j-1}^*,
    \end{align*}
    where $ \widetilde{C}_{j-1} = \frac{1}{N-1}d\widehat{V}_{j-1}^*  d\widehat{V}_{j-1} \in \R^{N \times N}$.
    Next, for fixed $ j \in \N $, we show that the eigenvectors of $ \widehat{C}_{j-1}  $ are also the eigenvectors of $ C_{j-1}  $.
    Indeed, if $ \phi \in \U $ satisfies $ \widehat{C}_{j-1}\phi = \lambda \phi $ with an eigenvalue $ \lambda \ge 0 $, we have
    \begin{align*}
        \widetilde{C}_{j-1}d\widehat{V}_{j-1}^* \phi 
        & = \frac{1}{N-1}d\widehat{V}_{j-1}^*  d\widehat{V}_{j-1} d\widehat{V}_{j-1}^*  = d\widehat{V}_{j-1}^*  \widehat{C}_{j-1}  \phi =   \lambda d\widehat{V}_{j-1}^* \phi.
    \end{align*}
    Hence, it follows that
    \begin{align*}
        C_{j-1} \phi 
        & = \frac{\alpha^2}{N-1}  d\widehat{V}_{j-1} (I_N + \alpha^2 \gamma^{-2}\widetilde{C}_{j-1})^{-1} d\widehat{V}_{j-1}^*  \phi = \frac{\alpha^2}{N-1}  d\widehat{V}_{j-1} \frac{1}{1 + \alpha^2\gamma^{-2}\lambda} d\widehat{V}_{j-1}^*  \phi \\
        & = \frac{\alpha^2}{1 + \alpha^2\gamma^{-2}\lambda} \widehat{C}_{j-1} \phi = \frac{\alpha^2 \lambda}{1 + \alpha^2\gamma^{-2}\lambda} \phi.
    \end{align*}
    Since the map $ \lambda \mapsto \frac{\alpha^2 \lambda}{1 + \alpha^2\gamma^{-2}\lambda} $ is monotonically increasing, we obtain the relation between the minimum eigenvalues
    \begin{align}
        \label{eq:lambda_min:update}
        \lambda^{min}_{j-1} = \frac{\alpha^2 \widehat{\lambda}^{min}_{j-1}}{1 + \frac{\alpha^2}{\gamma^2}\widehat{\lambda}^{min}_{j-1}}.
    \end{align}
    Combining \cref{ineq:lambda_min:chaos,eq:lambda_min:update}, we obtain the inequality for $ \widehat{\lambda}^{min}_j $.
    \begin{align*}
        \widehat{\lambda}^{min}_j \ge \frac{e^{-ah} \alpha^2 \widehat{\lambda}^{min}_{j-1}}{1 + \frac{\alpha^2}{\gamma^2}\widehat{\lambda}^{min}_{j-1}}.
    \end{align*}
    We now consider the following discrete dynamical system of the eigenvalue
    \begin{align*}
        \lambda_{n+1} = g(\lambda_n), \quad \lambda_0 > 0,
    \end{align*}
    where $ g(\lambda) = \frac{e^{-ah} \alpha^2 \lambda}{1 + \frac{\alpha^2}{\gamma^2}\lambda} $.
    Note that $ \lambda_n > 0 $ for $ n \in \N \cup \{0\} $.
    Let $ \lambda_\infty = \frac{\gamma^2}{\alpha^2}(e^{-ah}\alpha^2-1) $ be a fixed point of the dynamical system, i.e., $ \lambda_\infty = g(\lambda_\infty) $.
    Then, if $ e^{-ah} \alpha^2 > 1 $, the ratio $ \frac{g(\lambda)}{\lambda} $ satisfies $ \frac{g(\lambda)}{\lambda} \ge 1 $ (resp. $ < 1 $) for $ \lambda \le \lambda_\infty $ (resp. $ \lambda > \lambda_\infty $).
    Hence, we have $ \lim_{n \rightarrow \infty} \lambda_n = \lambda_\infty $.
    On the other hand, $ \lim_{n \rightarrow \infty} \lambda_n = 0 $ if $ e^{-ah}\alpha^2 \le 1 $.
    Therefore, we obtain the lower bound
    \begin{align}
        \label{ineq:lambda_min:lower}
        \widehat{\lambda}^{min}_j \ge \min\left\{\widehat{\lambda}^{min}_0, \frac{\gamma^2}{\alpha^2}(e^{-ah}\alpha^2-1)\right\} = \min\left\{ e^{-ah}\lambda_0, \frac{\gamma^2}{\alpha^2}(e^{-ah}\alpha^2-1)\right\} = \lambda_* > 0.
    \end{align}
    if and only if $ e^{-ah} \alpha^2 > 1 $.

    Finally, we establish the one-step inequality for $ \E[|e_j|^2] $.
    From \Cref{asmp:obs1}, equation \eqref{eq:mean_update_inflation} is reduced to
    \begin{align*}
        (I + \alpha^2 \gamma^{-2} \widehat{C}_j ) \overline{v}_j = \overline{\widehat{v}}_j + \alpha^2 \gamma^{-2}\widehat{C}_j  y_j.
    \end{align*}
    As in the proof of \Cref{thm:wellposed}, the error is devided into the two terms $ e_j = r_1 + r_2 $, where
    \begin{align}
        r_1 & = (I_\U + \alpha^2\gamma^{-2}\widehat{C}_j)^{-1}\widehat{e}_j,\label{eq:r1}\\
        r_2 & = (I_\U + \alpha^2\gamma^{-2}\widehat{C}_j)^{-1}\alpha^2\gamma^{-2}\widehat{C}_j  (y_j -  u_j),\label{eq:r2}
    \end{align}
    and $ \widehat{e}_j = \overline{\widehat{v}}_j - u_j $.
    From \cref{lem:i+a:succeq}, we have $ |(I_\U + \alpha^2\gamma^{-2}\widehat{C}_j )^{-1} \alpha^2\gamma^{-2}\widehat{C}_j | \le 1 $, and hence,
    \begin{align}
        \label{thm:error_bound:pf:r2}
        |r_2| \le |y_j -  u_j| = |\xi_j|.
    \end{align}
    From this lower bound of the minimum eigenvalue \eqref{ineq:lambda_min:lower}, we have $ |(I_\U + \alpha^2\gamma^{-2}\widehat{C}_j)^{-1}| \le (1+\frac{\alpha^2}{\gamma^2} \lambda_*)^{-1} $.
    Hence,
    \begin{align}
        \label{thm:error_bound:pf:Pm_r1}
        | r_1|^2 \le \frac{1}{(1+\frac{\alpha^2}{\gamma^2} \lambda_*)^2} | \widehat{e}_j|^2 \le  \frac{e^{2(\beta + \epsilon) h}}{(1+\frac{\alpha^2}{\gamma^2} \lambda_*)^2}(| e_{j-1}|^2 + D) - \frac{D}{(1+\frac{\alpha^2}{\gamma^2} \lambda_*)^2},
    \end{align}
    where \Cref{lem:chaos_proj} is used for any $ \epsilon > 0 $.

    As in the proof of \Cref{thm:wellposed}, we treat $ r_1 $ and $ r_2 $ separately when computing the expectation.
    Therefore, we obtain
    \begin{align*}
        \E[|e_j|^2]
        & = \E[|r_1|^2] + \E[|r_2|^2] \le \E\left[| r_1|^2\right] + \E[|\xi_j|^2] = \E\left[| r_1|^2 \right]+ m\gamma^2 \\
        & \le \frac{e^{2(\beta + \epsilon) h}}{(1+\frac{\alpha^2}{\gamma^2} \lambda_*)^2} (\E[| e_{j-1}|^2] + D) - \frac{D}{(1+\frac{\alpha^2}{\gamma^2} \lambda_*)^2} + m \gamma^2\\
        & \le \theta (\E[|e_{j-1}|^2] + D) - \frac{D}{(1+\frac{\alpha^2}{\gamma^2} \lambda_*)^2} + m \gamma^2 = \theta (\E[|e_{j-1}|^2] + D) + \mathcal{E},
    \end{align*}
    where $ \mathcal{E} = m \gamma^2 - \Theta D $.
    Then, we have
    \begin{align*}
        \E[|e_j|^2] + D \le \theta (\E[|e_{j-1}|^2] + D) + \mathcal{E} + D.
    \end{align*}
    Applying this inequality repeatedly, we have
    \begin{align*}
        \E[|e_j|^2] 
        & \le \theta^j (\E[|e_0|^2] + D) + (\mathcal{E} + D) \frac{1 - \theta^j}{1-\theta} - D \\
        & = \theta^j (\E[|e_0|^2] + D) + \left(m \gamma^2 + (1-\Theta)D\right) \frac{1 - \theta^j}{1-\theta} - D \\
        & = \theta^j (\E[|e_0|^2] + D) + m \gamma^2 \frac{1 - \theta^j}{1-\theta} + \left(\frac{(1 - \theta^j)(1 - \Theta)}{1 - \theta} - 1\right)D.
    \end{align*}
    Moreover, if $ \theta < 1 $, \eqref{ineq:main:result_limsup} holds in the limit of $ j \rightarrow \infty $.
\end{proof}

\begin{remark}
    The explicit condition about $ \alpha_0 $ so that $ \theta = (1+\frac{\alpha^2}{\gamma^2} \lambda_*)^{-2} e^{2(\beta + \epsilon) h} < 1 $ is given by
    \begin{align}
        \label{rem:ineq1}
        \alpha_0 = \max\{\lambda_0^{-\frac{1}{2}} \gamma e^{ah} (e^{(\beta + \epsilon)h} - 1)^{\frac{1}{2}}, e^{\frac{1}{2}(a + \beta + \epsilon)h}\}.
    \end{align}
    Indeed, the condition $ \theta < 1 $ is simplified to
    \begin{align*}
        e^{(\beta + \epsilon) h} - 1 < \frac{\alpha^2}{\gamma^2} \lambda_* = \min\left\{\frac{\alpha^2}{\gamma^2}e^{-ah}\lambda_0, e^{-ah}\alpha^2-1\right\}.
    \end{align*}
    Hence, \cref{rem:ineq1} is obtained by solving the inequalities for $ \alpha $.
\end{remark}

\begin{corollary}
    Under the same assumptions of \Cref{thm:error_bound}, we consider the accurate observation limit (i.e., $ \gamma \rightarrow 0 $).
    Then, the filtering error is the order of the observation noise.
    \begin{align*}
        \limsup_{j \rightarrow \infty} \E[|e_j|^2] = O(\gamma^2).
    \end{align*}
\end{corollary}
\begin{proof}
    It is clear that $ \Theta = \left(\frac{\gamma^2}{\gamma^2 + \alpha^2 \lambda_*}\right)^2 = O(\gamma^4) $ and $ \theta = O(\gamma^4) $.
    Then, we have
    \begin{align*}
        \frac{1-\Theta}{1-\theta} - 1 & = (1 - \Theta)(1 + \theta + O(\Theta^2)) - 1 = 1 - \Theta + \theta + O(\Theta^2) - 1 \\
        & = \Theta(e^{2(\beta + \epsilon)h} - 1) + O(\Theta^2) = O(\gamma^4).
    \end{align*}
    Therefore, it follows that
    \begin{align*}
        \limsup_{j \rightarrow \infty} \E[|e_j|^2] = \frac{m\gamma^2}{1 - \theta} + \left(\frac{1-\Theta}{1-\theta} - 1\right)D = m\gamma^2 (1 + O(\gamma^4)) + O(\gamma^4) = O(\gamma^2).
    \end{align*}
\end{proof}

\section{Summary and discussion}
\label{sec:discussion}
The theoretical aspects of the ETKF are investigated for the infinite-dimensional dynamics.
We obtain that the filtering error of the ETKF is bounded for any finite time without a covariance inflation.
In addition, an appropriate multiplicative covariance inflation to the ETKF ensures the uniform-in-time error bound on a finite-dimensional state space.
Furthermore, we determin the minimum value of the inflation parameter sufficient to obtain the uniform-in-time error bound.
As a corollary, the leading order of the error bound is equal to the scale $ \gamma $ of observation noises, indicating that the accuracy of the state estimation is effectively improved by the accurate observation.
These results are relevant to those for the PO method with the additive inflation \cite{kellyWellposednessAccuracyEnsemble2014b}.

Let us discuss the future directions.
Our analysis is limited to when the ensemble size $ N $ is larger than or equal to the state space dimension $ m $ to prevent a degenerated ensemble covariance.
Hence, we need to establish the error bound for the ETKF when $ N $ is smaller than $ m \le \infty $.
To this end, we can consider some detailed properties of the model dynamics in the present analysis.
For instance, it is useful to identify unstable directions of the filtering error where the error grows through the evolution of the model dynamics.
It is known that many dissipative dynamical systems including the two-dimensional Navier-Stokes equations have a finite number of such unstable directions \cite{constantinGlobalLyapunovExponents1985,constantinNavierstokesEquations1988,temamInfinitedimensionalDynamicalSystems1997}.
Hence, it is reasonable to assume that the dynamics have a finite number of unstable directions, and perturbations to the trajectory in the other directions decay through its evolution.
In addition, we assume the $ N $ is larger than or equal to the finite number of unstable directions.
If the ETKF appropriately reduces the error of these directions under these assumptions, the dynamics reduces the error of the other directions, and we obtain a better error bound.


\section*{Acknowledgments}
We used an AI tool to edit or polish the authors' written text for spelling, grammar, or general style.

\bibliographystyle{siamplain}
\bibliography{source/references}

\appendix
\section{Linear operators on Hilbert space}
Let $ \H $ be a Hilbert space and $ I = I_\H $.
\begin{lemma}
    \label{lem:i+a}
    Let $ A \in \L(\H) $.
    If $ I + A $ is invertible, then we have
    \begin{align}
        \label{lem:i+a:inv}
        (I + A)^{-1} = I - (I+A)^{-1} A.
    \end{align}
    Especially, if $ A \in \L_{sa}(\H) $ and $ A \succ 0 $,
    \begin{align}
        \label{lem:i+a:succeq}
        0 \preceq A(A + I)^{-1} = (A + I)^{-1}A \preceq I, \quad 0 \preceq (A+I)^{-1} \preceq I.
    \end{align}
\end{lemma}
\begin{proof}
    \cref{lem:i+a:inv} is easily confirmed by
    \begin{align*}
        LHS = (I+A)^{-1} (I + A - A) = I - (I+A)^{-1} A.
    \end{align*}
    Next, for \cref{lem:i+a:succeq}, \cref{lem:i+a:inv} yields $ A(A + I)^{-1} = (A + I)^{-1}A $.
    The inequalities hold from the spectral mapping theorem.
\end{proof}

\begin{lemma}
    \label{lem:invertible_operator}
    Let $ \Gamma: \H \rightarrow \H $ be invertible and $ V \in \H^N $.
    Then, the operator $ VV^* + \Gamma $ is invertible, and
    \begin{align}
        \label{lem:invertible_operator:pre}
        (I + V^*\Gamma^{-1} V)^{-1}V^*\Gamma^{-1} = V^*(VV^* + \Gamma)^{-1}.
    \end{align}
    Furthermore, 
    \begin{align}
        \label{lem:invertible_operator:main}
        (I + V^*\Gamma^{-1} V)^{-1} = I - V^*(VV^* + \Gamma)^{-1}V.
    \end{align}
\end{lemma}
\begin{proof}
    $ VV^* + \Gamma $ is invertible owing to $ VV^* \succeq 0 $ and $ \Gamma \succ 0 $.
    Then, we have
    \begin{align*}
        V^* \Gamma^{-1} (VV^* + \Gamma) = V^*\Gamma^{-1} VV^* + V^* = (I + V^*\Gamma^{-1}V)V^*.
    \end{align*}
    This is equivalent to \cref{lem:invertible_operator:pre}.
    For \cref{lem:invertible_operator:main}, the assertion \cref{lem:invertible_operator:pre} yields
    \begin{align*}
        (I + V^*\Gamma^{-1} V)^{-1}V^*\Gamma^{-1}V = V^*(VV^* + \Gamma)^{-1}V.
    \end{align*}
    Hence, we have
    \begin{align*}
        I - V^*(VV^* + \Gamma)^{-1}V = I - (I + V^*\Gamma^{-1} V)^{-1}V^*\Gamma^{-1}V = (I + V^*\Gamma^{-1} V)^{-1},
    \end{align*}
    where the last equality holds from \cref{lem:i+a:inv}.
\end{proof}

\end{document}